\documentclass[11pt]{article}
\usepackage{pifont}
\usepackage{amssymb}
\usepackage{amsfonts}
\usepackage{amscd}
\usepackage{epsfig}
\usepackage{amsmath,amstext,amsthm,amsbsy,amssymb}

\newcommand{\bc}{{\mathbb C}}
\newcommand{\br}{{\mathbb R}}
\newcommand{\bh}{{\mathbb H}}

\newcommand{\mtx}[4]

\newtheorem{thm}{Theorem}[section]

\newtheorem{lem}{Lemma}[section]
\newtheorem{prop}{Proposition}[section]

\newtheorem{dfn}{Definition}

\begin{document}

\title{J{\o}rgensen's  inequality for quaternionic  hyperbolic n-space}
\author{Wensheng Cao\thanks{Supported by NSFs of China (No.10801107, No.10671004 ) and
NSF of Guangdong Province (No.8452902001000043)}\\
Department of Mathematics and Physics, Wuyi University,\\
Jiangmen, Guangdong 529020, P.R. China\\
e-mail:  wenscao@yahoo.com.cn}
\date{}
\maketitle

\bigskip
{\bf Abstract}: J{\o}rgensen's inequality gives a necessary
condition for a non-elementary two generator group of isometries of
real hyperbolic 2-space to be discrete. We give analogues of
J{\o}rgensen's inequality for non-elementary groups of isometries of
quaternionic hyperbolic n-space generated by two elements, one of
which is loxodromic.

\medskip

{\bf Mathematics Subject Classifications (2000):} 20H10, 22E40,
57S30.
\medskip

{\bf Keywords:}  quaternionic hyperbolic n-space,  J{\o}rgensen's
inequality.

\section{Introduction}

J{\o}rgensen's  inequality \cite{jor} gives a necessary condition
for a non-elementary two generator subgroup of $PSL(2,\bc)$ to be
discrete.   As a quantitative version of Margulis' lemma, this
inequality has been generalized by many people in different
contexts. Viewing  $PSL(2,\br)$  as the isometry group of complex
hyperbolic 1-space, $H_{\bc}^1$,  we can seek to generalize
J{\o}rgensen's inequality to higher dimensional complex hyperbolic
isometries.  Kamiya \cite{kam83,kam91} and Parker \cite{par92,par97}
gave generalizations of J{\o}rgensen's inequality to the two
generator subgroup of $PU(n,1)$ when one generator is Heisenberg
translation. By using the stable basin theorem, Basmajian and Miner
\cite{bm} generalized the J{\o}rgensen's  inequality to the two
generator subgroup of $PU(2, 1)$ when the generators are loxodromic
or boundary elliptic, and several other inequalities are due to
Jiang, Kamiya and Parker \cite{jkp} by using matrix method other
than the purely geometric method.   Jiang \cite{jp} and Kamiya
\cite{kp} gave generalizations of J{\o}rgensen's inequality to the
two generator subgroup of $PU(2,1)$ when one generator is Heisenberg
screw motion. The generalization also was done  in \cite{xj} for the
case when one generator is a regular elliptic element.

Following the research on complex hyperbolic space, I. Kim and J.
Parker opened up the study of quaternionic hyperbolic space in
\cite{kimp}.  They proved some basic facts about discreteness of two
generator subgroups, minimal volume of cusped quaternionic
manifolds, and laid some basic tools to study quaternionic
hyperbolic space.

It is naturally asked that theorems in complex hyperbolic space
can be generalized to quaternionic hyperbolic space.  There has
been several investigation  in this area. Kim \cite{kim} found
analogues in quaternionic hyperbolic 2-space of results in
\cite{jkp,jp}.  Cao \cite{caog} obtained analogue of
J{\o}rgensen's inequality for non-elementary groups of isometries
of quaternionic hyperbolic n-space generated by two elements, one
of which is elliptic.  Markham and Parker \cite{mp} studied
related problems in metric space.

 In this paper, we will give analogues of J{\o}rgensen's
inequality for non-elementary groups of isometries of quaternionic
hyperbolic n-space generated by two elements, one of which is
loxodromic.

 Our main result is

\begin{thm}\label{mainthm}
Let $g$ and $h$ be elements of $Sp(n,1)$ such that  $g$ is
loxodromic element with fixed points $u,v$. If either
\begin{equation}\label{cond} M_g(1+|[h(u),v,u,h(v)]|^{1/2})<1\ \ \mbox{or}\ \ M_g(1+|[h(u),u,v, h(v)]|^{1/2})<1, \end{equation}
then the group $\langle g,h \rangle$ is either elementary or not
discrete.
\end{thm}

We arrange this paper as follows. Section 2 contains some necessary
background materials of quaternionic hyperbolic geometry. Section 3
contains the definition of quaternionic cross-ratio and a lemma
which is crucial to prove our main result.  Section 4 aims at the
proof of Theorem \ref{mainthm}.

\section{Preliminaries}

Let $\bh$ denote the division ring of real quaternions. Elements of
$\bh$ have the form $z=z_1+z_2{\bf i}+z_3{\bf j}+z_4{\bf k}\in \bh$
where $z_i\in \br$ and
$$
{\bf i}^2 = {\bf j}^2 = {\bf k}^2 = {\bf i}{\bf j}{\bf k} = -1.
$$
Let $\overline{z}=z_1-z_2{\bf i}-z_3{\bf j}-z_4{\bf k}$ be the {\sl
conjugate} of $z$, and
$$
|z|= \sqrt{\overline{z}z}=\sqrt{z_1^2+z_2^2+z_3^2+z_4^2}
$$
be the {\sl modulus} of $z$. We define $\Re(z)=(z+\overline{z})/2$
to be the {\sl real part} of $z$, and $\Im(z)=(z-\overline{z})/2$ to
be the {\sl imaginary part} of $z$. Also
$z^{-1}=\overline{z}|z|^{-2}$ is the {\sl inverse} of $z$. Observe
that $\Re(wzw^{-1})=\Re(z)$ and $|wzw^{-1}|=|z|$ for all $z$ and $w$
in $\bh$. Two quaternions $z$ and $w$ are {\sl similar} if there
exists nonzero $q \in \bh$  such that $z=q w q^{-1}$. The similarity
class of $z$ is the set $\bigl\{ q z q^{-1} : q \in \bh-\{ 0\}
\bigr\}$.

We collect the following simple properties of quaternions as a
proposition for later use;  More properties of quaternions and
matrices of quaternions can be found in \cite{zhang}.

\begin{prop}\label{proquat}
For any two quaternions $a,b$, we have
$$\Re(ab)=\Re(ba)$$
and
$$2\Re(a)\Re(b)-\Re(ab)\leq |ab|.$$
\end{prop}

{\bf Proof.} Let $$a=a_1+a_2{\bf i}+a_3{\bf j}+a_4{\bf k},\ \
b=b_1+b_2{\bf i}+b_3{\bf j}+b_4{\bf k},$$ where $a_i,b_i\in \br,
i=1,\cdots,4.$ Then
$$\Re(ab)=a_1b_1-a_2b_2-a_3b_3-a_4b_4=\Re(ba)$$ and
$$2\Re(a)\Re(b)-\Re(ab)=a_1b_1+a_2b_2+a_3b_3+a_4b_4\leq
(a_1^2+a_2^2+a_3^2+a_4^2)^{1/2}(b_1^2+b_2^2+b_3^2+b_4^2)^{1/2}=|ab|.$$

In what follows, we give some necessary background materials of
quaternionic hyperbolic geometry. More details can be found in
\cite{chen,gold,kimp}.

Let $\bh^{n,1}$ be the quaternionic vector space of  quaternionic
dimension $n+1$  (so real dimension $4n+4$) with the quaternionic
Hermitian form
$$
\langle{\bf z},\,{\bf w}\rangle={\bf w}^*J{\bf z}=
\overline{w_1}z_1+\cdots+\overline{w_{n-1}}z_{n-1}-(\overline{w_{n}}z_{n+1}+\overline{w_{n+1}}z_{n}),
$$
where ${\bf z}$ and ${\bf w}$ are the column vectors in $\bh^{n,1}$
with entries $z_1,\cdots,z_{n+1}$ and $w_1,\cdots,w_{n+1}$
respectively, $\cdot^*$ denotes quaternionic Hermitian transpose and
$J$ is the Hermitian matrix
$$J=\left(
                  \begin{array}{ccc}
                    I_{n-1} & 0 & 0 \\
                    0 & 0 & -1 \\
                    0 & -1 & 0\\
                  \end{array}
                \right).$$
We define a {\sl unitary transformation} $g$ to be an automorphism
$\bh^{n,1}$, that is, a linear bijection such that $\langle g({\bf
z}),\,g({\bf w})\rangle=\langle{\bf z},\,{\bf w}\rangle$ for all
${\bf z}$ and ${\bf w}$ in $\bh^{n,1}$. We denote the group of all
unitary transformations by $Sp(n,1)$.

Following Section 2 of \cite{chen}, let
\begin{eqnarray*}
V_0 & = & \Bigl\{{\bf z} \in  \bh^{n,1}-\{0\}:
\langle{\bf z},\,{\bf z}\rangle=0\Bigr\} \\
V_{-} &  = & \Bigl\{{\bf z} \in \bh^{n,1}:\langle{\bf z},\,{\bf
z}\rangle<0\Bigr\}.
\end{eqnarray*}
It is obvious that $V_0$ and $V_{-}$ are invariant under $Sp(n,1)$.
We define $V^s$ to be $V^s=V_{-}\cup  V_0$. Let $P:V^s\to
P(V^s)\subset \bh^{n}$ be the right  projection map defined by
$$
P(z_1,\cdots,z_n,  z_{n+1})^t=(z_1z_{n+1}^{-1},\cdots,z_n
z_{n+1}^{-1})^t,
$$
where $\cdot^t$ denotes the  transpose.

 We define $H_{\bh}^n=P(V_-)$
and $\partial H_{\bh}^n=P(V_0)$.
  We call $H_{\bh}^n$ the siegel domain  model of quaternionic hyperbolic $n$-space.  The Bergman metric on $H_{\bh}^n$ is given by the distance formula
 $$\cosh^2\frac{\rho(z,w)}{2}=\frac{\langle{\bf z},\,{\bf
w}\rangle \langle{\bf w},\,{\bf z}\rangle}{\langle{\bf z},\,{\bf
z}\rangle \langle{\bf w},\,{\bf w}\rangle},\ \ \mbox{where}\ \ z,w
\in H_{\bh}^n, \ \  {\bf z}\in P^{-1}(z),{\bf w}\in P^{-1}(w).$$
 The  isometry group of $H_{\bh}^n$ with respect to the
Bergman metric is the projective unitary group $PSp(n,1)$ and acts
on $P(\bh^{n,1})$ by matrix multiplication.  Here, we adopt the
convention that the action of $Sp(n,1)$ on $H_{\bh}^n$ is left
action and the action of projectivization of $Sp(n,1)$ is right
action.

If $g\in Sp(n,1)$,  by definition, $g$ preserves the Hermitian form.
Hence
$$
{\bf w}^*J{\bf z}=\langle{\bf z},\,{\bf w}\rangle= \langle g{\bf
z},\,g{\bf w}\rangle ={\bf w}^* g^*J g{\bf z}
$$
for all ${\bf z}$ and ${\bf w}$ in $\bh^{n,1}$. Letting ${\bf z}$
and ${\bf w}$ vary over a basis for $\bh^{n,1}$, we see that $J=
g^*J g$. From this we find $ g^{-1}=J^{-1} g^*J$. That is:
\begin{equation}g^{-1}=\left(
  \begin{array}{ccc}
     A^*& -\theta^*& -\gamma^* \\
    -\beta^* & \overline{a_{n+1,n+1}}& \overline{a_{n,n+1}}\\
    -\alpha^* & \overline{a_{n+1,n}}& \overline{a_{n,n}}\\
    \end{array}
\right) \end{equation}

if

\begin{equation}
 g=(a_{ij})_{i,j=1,\cdots,n+1}=\left(
  \begin{array}{ccc}
     A& \alpha& \beta \\
    \gamma & a_{n,n}& a_{n,n+1}\\
    \theta & a_{n+1,n}& a_{n+1,n+1}\\
    \end{array}
\right)\in Sp(n,1).\end{equation}

 Using the identities $gg^{-1}=g^{-1}g=I_{n+1}$ we obtain:

\begin{eqnarray}
 \label{AA}&&AA^*-\alpha \beta^*-\beta \alpha^*=I_{n-1},\\
 && -A\theta^*+\alpha \overline{a_{n+1,n+1}}+\beta \overline{a_{n+1,n}}= 0,\\
 &&-A\gamma^*+\alpha \overline{a_{n,n+1}}+\beta \overline{a_{n,n}}= 0,\\
 &&-\gamma\theta^*+a_{n,n}\overline{a_{n+1,n+1}}+a_{n,n+1}\overline{a_{n+1,n}}= 1,\\
&&\label{gg}-\gamma\gamma^*+a_{n,n}\overline{a_{n,n+1}}+a_{n,n+1}\overline{a_{n,n}}=0,\\
&&\label{tt}-\theta\theta^*+a_{n+1,n}\overline{a_{n+1,n+1}}+a_{n+1,n+1}\overline{a_{n+1,n}}=0,\\
&&A^*A-\theta^* \gamma-\gamma^* \theta=I_{n-1},\\
&& A^*\alpha-\theta^*a_{n,n}-\gamma^*a_{n+1,n}= 0,\\
&&A^*\beta-\theta^*a_{n,n+1}-\gamma^* a_{n+1,n+1}= 0,\\
&&\label{bba}-\beta^*\alpha+\overline{a_{n+1,n+1}}a_{n,n}+\overline{a_{n,n+1}}a_{n+1,n}= 1,\\
&&\label{bb} -\beta^*\beta+\overline{a_{n+1,n+1}}a_{n,n+1}+\overline{a_{n,n+1}}a_{n+1,n+1}=0,\\
&&\label{aa}-\alpha^*\alpha+\overline{a_{n+1,n}}a_{n,n}+\overline{a_{n,n}}a_{n+1,n}=0.
\end{eqnarray}

For a non-trivial element $g$ of $Sp(n,1)$, we say that $g$ is
\begin{itemize}
\item[(i)] {\sl elliptic} if it has a fixed point in  $H_{\bh}^n$;
\item[(ii)] {\sl parabolic} if it has exactly one fixed point which lies in
$\partial H_{\bh}^n$;
\item[(iii)] {\sl loxodromic} if it has exactly two fixed points which
lie in $\partial H_{\bh}^n$.
\end{itemize}

 A  subgroup $G$ of
$Sp(n,1)$ is called {\it non-elementary} if it contains two
non-elliptic elements of infinite order with distinct fixed points;
Otherwise $G$ is called {\it elementary}.

As in complex hyperbolic n-space, we have the following
proposition classifying  elementary subgroups of $Sp(n, 1)$.
\begin{prop}(cf.\cite[Lemma 2.4]{caobull} \label{element}
(i)\quad If  $G$ contains a parabolic element but no loxodromic
element, then $G$ is elementary if and only if it fixes a point in
$\partial H_{\bh}^n$;

(ii)\quad If  $G$ contains a loxodromic element, then $G$ is
elementary if and only if it fixes a point in $\partial H_{\bc}^n$
or a point-pair $\{x, y\}\subset \partial H_{\bh}^n$;

(iii)\quad $G$ is purely elliptic, i.e., each non-trivial element
of $G$ is elliptic, then  $G$ is elementary and fixes a point in
$\overline{H_{\bh}^n}$.
\end{prop}

Let  $q_0, q_\infty\in \partial H_{\bh}^n$ stand for the images of
$(0,\cdots,0,1), (0,\cdots,0,1,0)\in \bh^{n,1}$ under  the
projection map $P$, respectively, and $$G_0=\{g\in Sp(n,1):
g(q_0)=q_0\}, G_{\infty}=\{g\in Sp(n,1): g(q_\infty)=q_\infty\},
G_{0,\infty}=G_0\cap G_{\infty}.$$ Then we have the following
three propositions, see \cite[Lemma 3.3.1]{chen}.

\begin{prop}
If $g\in G_{\infty}$, then $g$ is of form
 \begin{equation}\label{gw}\left(
  \begin{array}{ccc}
    A & 0 & a \\
    b & \lambda & s \\
    0 & 0 & \mu \\
  \end{array}
\right), \end{equation}
 where $\mu, \lambda, s\in \bh$, ${\overline \mu
}\lambda
=1,\;\Re(\overline{\mu}s)=\frac{1}{2}\left|a\right|^2,\;b=\lambda
a^* A$ and $A\in U(n-1;\bh).$
\end{prop}

\begin{prop}
If $g\in G_{0}$, then $g$ is of form
 \begin{equation}\label{g0}\left(
  \begin{array}{ccc}
    A & a & 0 \\
    0 & \mu & 0 \\
    b & s & \lambda \\
  \end{array}
\right),\end{equation} where $\mu, \lambda, s\in \bh$, ${\overline
\mu }\lambda
=1,\;\Re(\overline{\mu}s)=\frac{1}{2}\left|a\right|^2,\;b=\lambda
a^* A$ and $A\in U(n-1;\bh).$
\end{prop}

\begin{prop}
If $g\in G_{0,\infty}$, then $g$ is of form
 \begin{equation}\label{g0w}\left(
  \begin{array}{ccc}
    A & 0 & 0 \\
    0 & \lambda & 0 \\
    0 & 0 & \mu \\
  \end{array}
\right),\end{equation} where $\mu, \lambda, \in \bh$, ${\overline
\mu }\lambda =1$ and $A\in U(n-1;\bh).$
\end{prop}

If $g\in Sp(n,1)$ is  loxodromic, then $g$  conjugates to
\begin{equation}\label{loxo}
diag(\lambda_1,\  \lambda_2,\  \cdots,\  \lambda_{n+1})\in
G_{0,\infty},
\end{equation}
 where $\lambda_i, \ i=1,\cdots,n-1$  are right eigenvalues of quaternionic matrix $g$  with norm 1 and
$|\lambda_n|\neq 1, \overline{\lambda_{n}}\lambda_{n+1}=1.$

Right eigenvalues are conjugacy invariants for matrices $g\in
Sp(n,1)$, see \cite{cpw,caope, zhang}. In fact, if $g{\bf z}={\bf
z}t$ then
$$
(hgh^{-1})( h{\bf z}) = h g({\bf z})=( h{\bf z})t.
$$
If $t$ is a right eigenvalue of $g$ then so is any quaternion in the
similarity class of $t$. In order to see this, observe that if $
g{\bf z}={\bf z}t$ and $q$ is any non-zero quaternion then
$$
g({\bf z}q)= g{\bf z}q={\bf z}tq=({\bf z}q)(q^{-1}tq).
$$

For any loxodromic element $g\in Sp(n,1)$,  after conjugating to the
form (\ref{loxo}) we define
\begin{equation}\delta(g)=\max_{1\leq i \leq n-1 }|\lambda_i-1|\end{equation}
and

\begin{equation}
M_g=2\delta(g)+|\lambda_n-1|+|\lambda_{n+1}-1|.
\end{equation}
Then we have the following lemma.
\begin{lem}
If $g\in Sp(n,1)$  is  loxodromic, then $\delta(g)$  and $M_g$ are
conjugacy invariants.
\end{lem}

\section{Main Lemma}

Cross-ratios were generalized to complex hyperbolic space by
Kor\'anyi and Reimann \cite{km}. We will generalize this definition
of complex cross-ratio to the non commutative quaternion ring.

\begin{dfn}  The quaternionic
cross-ratio of four points  $z_1, z_2, w_1, w_2$ in
$\overline{H_{\bh}^n}$  is defined as :
\begin{equation}
 [z_1, z_2, w_1, w_2]=
\langle{\bf w_1},\,{\bf z_1}\rangle \langle{\bf w_1},\,{\bf
z_2}\rangle^{-1} \langle{\bf w_2},\,{\bf z_2}\rangle \langle{\bf
w_2},\,{\bf z_1}\rangle^{-1},
\end{equation}
where ${\bf z}_i=\in P^{-1}(z_i)$ and ${\bf w_i}\in P^{-1}(w_i)$ for
$i=1,\,2$.
\end{dfn}

Obviously, The quaternionic cross-ratio $[z_1, z_2, w_1, w_2]$
depends on the choice of ${\bf z}_1\in P^{-1}(z_1)$. However, its
absolute value
\begin{equation}
 \bigl|[z_1, z_2, w_1, w_2]\bigr|
=\frac{|\langle{\bf w_1},\,{\bf z_1}\rangle\langle{\bf w_2},\,{\bf
z_2}\rangle|} {| \langle{\bf w_1},\,{\bf z_2}\rangle\langle{\bf
w_2},\,{\bf z_1}\rangle|}
\end{equation}
is independent of the preimage of $z_i$ and $w_i$ in $\bh^{n,1}$.

 Let \begin{equation}\label{elementh}
 h=(a_{ij})_{i,j=1,\cdots,n+1}=\left(
  \begin{array}{ccc}
     A& \alpha& \beta \\
    \gamma & a_{n,n}& a_{n,n+1}\\
    \theta & a_{n+1,n}& a_{n+1,n+1}\\
    \end{array}
\right)\in Sp(n,1).\end{equation} Then

 \begin{equation}\label{ratio1}
 |[h(q_{\infty}),q_{0},q_{\infty}, h(q_{0})]|=|a_{n+1,n}a_{n,n+1}|
 \end{equation}
and
 \begin{equation}\label{ratio2}
 |[h(q_{\infty}),q_{\infty},q_{0},
 h(q_{0})]|=|a_{n,n}a_{n+1,n+1}|.
 \end{equation}

In order to find the discrete condition, we will use the absolute
value of this quaternionic cross-ratio.

The following lemma is crucial for us.

\begin{lem}\label{lemcu}
 Let $h$ be of the form (\ref{elementh}). Then

\begin{equation}\label{ba}
|\beta^*\alpha|\leq
2|a_{n,n}a_{n+1,n+1}|^{1/2}|a_{n,n+1}a_{n+1,n}|^{1/2},
\end{equation}

\begin{equation}\label{gt}
|\gamma\theta^*|\leq
2|a_{n,n}a_{n+1,n+1}|^{1/2}|a_{n,n+1}a_{n+1,n}|^{1/2},
\end{equation}

\begin{equation}\label{nn}
|a_{n,n}a_{n+1,n+1}|^{1/2}\leq |a_{n,n+1}a_{n+1,n}|^{1/2}+1,
\end{equation}

\begin{equation}\label{nd} |a_{n,n+1}a_{n+1,n}|^{1/2}\leq
|a_{n,n}a_{n+1,n+1}|^{1/2}+1,
\end{equation}

\begin{equation}\label{cl1}
|a_{n,n}a_{n+1,n+1}|^{1/2}+|a_{n,n+1}a_{n+1,n}|^{1/2}\geq 1.
\end{equation}
\end{lem}

 {\bf Proof.} By (\ref{bb}) and (\ref{aa}), we have
\begin{eqnarray}\label{bbaa}|\beta^*\alpha|^2&\leq& |\beta^*\beta||\alpha^*\alpha|
=2\Re(\overline{a_{n+1,n+1}}a_{n,n+1})2\Re(\overline{a_{n+1,n}}a_{n,n})\\
&\leq& 4|a_{n,n}a_{n+1,n+1}||a_{n,n+1}a_{n+1,n}|.
\end{eqnarray}

By (\ref{gg}) and (\ref{tt}), we have
\begin{eqnarray*}|\gamma\theta^*|^2&\leq& |\gamma\gamma^*|  |\theta\theta^*|=2\Re(a_{n,n}\overline{a_{n,n+1}})2\Re(a_{n+1,n}\overline{a_{n+1,n+1}})\\
&\leq& 4|a_{n,n}a_{n+1,n+1}||a_{n,n+1}a_{n+1,n}|.
\end{eqnarray*}

By (\ref{bba}) and  (\ref{bbaa}), we have that
\begin{eqnarray*}&&(1-\overline{a_{n+1,n+1}}a_{n,n}-\overline{a_{n,n+1}}a_{n+1,n})(1-\overline{a_{n,n}}a_{n+1,n+1}-\overline{a_{n+1,n}}a_{n,n+1})\\
&&=1+|a_{n+1,n+1}a_{n,n}|^2+|a_{n,n+1}a_{n+1,n}|^2+2\Re(\overline{a_{n+1,n+1}}a_{n,n}\overline{a_{n+1,n}}a_{n,n+1})\\
&&-2\Re(\overline{a_{n,n}}a_{n+1,n+1})-2\Re(\overline{a_{n+1,n}}a_{n,n+1})\\
&&=|\beta^*\alpha|^2\leq (\beta^*\beta)
(\alpha^*\alpha)=4\Re(\overline{a_{n+1,n+1}}a_{n,n+1})\Re(\overline{a_{n+1,n}}a_{n,n}).
\end{eqnarray*}
Therefore
\begin{eqnarray*}&&2\Re(\overline{a_{n,n}}a_{n+1,n+1})+2\Re(\overline{a_{n+1,n}}a_{n,n+1})\geq 1+|a_{n+1,n+1}a_{n,n}|^2+|a_{n,n+1}a_{n+1,n}|^2\\
&&\hspace{4mm}+2\Re(\overline{a_{n+1,n+1}}a_{n,n}\overline{a_{n+1,n}}a_{n,n+1})-4\Re(\overline{a_{n+1,n+1}}a_{n,n+1})\Re(\overline{a_{n+1,n}}a_{n,n})\\
&&=1+|a_{n+1,n+1}a_{n,n}|^2+|a_{n,n+1}a_{n+1,n}|^2\\
&&\hspace{4mm} -2(2\Re(a_{n,n+1}\overline{a_{n+1,n+1}})\Re(a_{n,n}\overline{a_{n+1,n}})-\Re(a_{n,n+1}\overline{a_{n+1,n+1}}a_{n,n}\overline{a_{n+1,n}}))\\
&&\geq
1+|a_{n+1,n+1}a_{n,n}|^2+|a_{n,n+1}a_{n+1,n}|^2-2|a_{n+1,n+1}a_{n,n}||a_{n,n+1}a_{n+1,n}|\\
&&=1+(|a_{n+1,n+1}a_{n,n}|-|a_{n,n+1}a_{n+1,n}|)^2.
\end{eqnarray*}
The last inequality follows from Proposition \ref{proquat}. Hence
$$ 2|a_{n,n}a_{n+1,n+1}|+2|a_{n+1,n}a_{n,n+1}|\geq
1+(|a_{n+1,n+1}a_{n,n}|-|a_{n,n+1}a_{n+1,n}|)^2.$$ That is

$$ 4|a_{n+1,n}a_{n,n+1}|\geq
(|a_{n+1,n+1}a_{n,n}|-|a_{n,n+1}a_{n+1,n}|-1)^2.$$ Taking the square
root gives
$$ -2|a_{n+1,n}a_{n,n+1}|^{1/2}\leq
|a_{n+1,n+1}a_{n,n}|-|a_{n,n+1}a_{n+1,n}|-1\leq
2|a_{n+1,n}a_{n,n+1}|^{1/2}.$$ Thus
$$ (1-|a_{n+1,n}a_{n,n+1}|^{1/2})^2 \leq
|a_{n+1,n+1}a_{n,n}|\leq  (1+|a_{n+1,n}a_{n,n+1}|^{1/2})^2 .$$
Taking the square root  in the above formula gives the desired
inequalities.

The proof is complete.

\section{The proof of Theorem \ref{mainthm}}

{\bf The proof of Theorem \ref{mainthm}.} \ \    Since (\ref{cond})
is invariant under conjugation, we may assume that $g$ is of  form
(\ref{loxo}) and $h$ is  of the form (\ref{elementh}).
 By (\ref{ratio1}) and (\ref{ratio2}), conditions (\ref{cond}) can be rewritten as
 \begin{equation}\label{rcond}  M_g(1+|a_{n,n+1}a_{n+1,n}|^{1/2})<1\ \ \mbox{or}\
 \ \
 M_g(1+|a_{n+1,n+1}a_{n,n}|^{1/2})<1.
 \end{equation}
Let $h_0 = h$  and  $h_{k+1} = h_kgh_k^{-1}$.  We write
$$h_k=(a^{(k)}_{ij})_{i,j=1,\cdots,n+1}=\left(
  \begin{array}{ccc}
     A^{(k)}& \alpha^{(k)} & \beta^{(k)} \\
    \gamma^{(k)} & a^{(k)}_{n,n}& a^{(k)}_{n,n+1}\\
    \theta^{(k)} & a^{(k)}_{n+1,n}& a^{(k)}_{n+1,n+1}\\
    \end{array}
\right).$$
 Then
\begin{eqnarray*}
h_{k+1}&=& \left(
  \begin{array}{ccc}
     A^{(k+1)}& \alpha^{(k+1)} & \beta^{(k+1)} \\
    \gamma^{(k+1)} & a^{(k+1)}_{n,n}& a^{(k+1)}_{n,n+1}\\
    \theta^{(k+1)} & a^{(k+1)}_{n+1,n}& a^{(k+1)}_{n+1,n+1}\\
    \end{array}
\right)\\
 &=& \left(
  \begin{array}{ccc}
     A^{(k)}& \alpha^{(k)} & \beta^{(k)} \\
    \gamma^{(k)} & a^{(k)}_{n,n}& a^{(k)}_{n,n+1}\\
    \theta^{(k)} & a^{(k)}_{n+1,n}& a^{(k)}_{n+1,n+1}\\
    \end{array}
\right)\left(
       \begin{array}{ccc}
         L & 0 & 0\\
         0 & \lambda_{n}& 0 \\
         0& 0 & \lambda_{n+1}\\
             \end{array}
     \right)\left(
  \begin{array}{ccc}
     (A^{(k)})^*& -(\theta^{(k)})^*& -(\gamma^{(k)})^* \\
    -(\beta^{(k)})^* & \overline{a^{(k)}_{n+1,n+1}}& \overline{a^{(k)}_{n,n+1}}\\
    -(\alpha^{(k)})^* & \overline{a^{(k)}_{n+1,n}}& \overline{a^{(k)}_{n,n}}\\
    \end{array}
\right),
\end{eqnarray*}
where $L=diag(\lambda_1,\  \lambda_2,\  \cdots,\  \lambda_{n-1}).$
Therefore
 \begin{eqnarray}\label{tolam1} &&a^{(k+1)}_{n,n}
 =-\gamma^{(k)}L(\theta^{(k)})^*+
a^{(k)}_{n,n}\lambda_n\overline{a^{(k)}_{n+1,n+1}}+a^{(k)}_{n,n+1}\lambda_{n+1}\overline{a^{(k)}_{n+1,n}},\\
&&\label{nn+1}a^{(k+1)}_{n,n+1}=-\gamma^{(k)}L(\gamma^{(k)})^*+
a^{(k)}_{n,n}\lambda_n\overline{a^{(k)}_{n,n+1}}+a^{(k)}_{n,n+1}\lambda_{n+1}\overline{a^{(k)}_{n,n}},\\
&&\label{n+1n}a^{(k+1)}_{n+1,n}=-\theta^{(k)}L(\theta^{(k)})^*+a^{(k)}_{n+1,n}\lambda_n\overline{a^{(k)}_{n+1,n+1}}
+a^{(k)}_{n+1,n+1}\lambda_{n+1}\overline{a^{(k)}_{n+1,n}},\\
\label{tolam2}&&a^{(k+1)}_{n+1,n+1}=
-\theta^{(k)}L(\gamma^{(k)})^*+a^{(k)}_{n+1,n}\lambda_n\overline{a^{(k)}_{n,n+1}}
+a^{(k)}_{n+1,n+1}\lambda_{n+1}\overline{a^{(k)}_{n,n}}.
\end{eqnarray}

By (\ref{gg}) and (\ref{nn+1}), we have
 \begin{eqnarray*}
|a^{(k+1)}_{n,n+1}|&=&|\gamma^{(k)}(I_{n-1}-L)(\gamma^{(k)})^*+
a^{(k)}_{n,n}(\lambda_n-1)\overline{a^{(k)}_{n,n+1}}+a^{(k)}_{n,n+1}(\lambda_{n+1}-1)\overline{a^{(k)}_{n,n}}|\\
&=&|\sum_{j=1}^{n-1}a^{(k)}_{n,j}(1-\lambda_{j})\overline{a^{(k)}_{n,j}}
+a^{(k)}_{n,n}(\lambda_n-1)\overline{a^{(k)}_{n,n+1}}+a^{(k)}_{n,n+1}(\lambda_{n+1}-1)\overline{a^{(k)}_{n,n}}|\\
&\leq
&\delta(g)\sum_{j=1}^{n-1}|a^{(k)}_{n,j}|^2+(|\lambda_n-1|+|\lambda_{n+1}-1|)|a^{(k)}_{n,n+1}a^{(k)}_{n,n}|\\
&=&
\delta(g)2\Re(a^{(k)}_{n,n}\overline{a^{(k)}_{n,n+1}})+(|\lambda_n-1|+|\lambda_{n+1}-1|)|a^{(k)}_{n,n+1}a^{(k)}_{n,n}|\\
&\leq
&(2\delta(g)+|\lambda_n-1|+|\lambda_{n+1}-1|)|a^{(k)}_{n,n+1}a^{(k)}_{n,n}|\\
&=& M_g|a^{(k)}_{n,n+1}a^{(k)}_{n,n}|.
\end{eqnarray*}

Similarly, by (\ref{tt})  and (\ref{n+1n})  we have
\begin{eqnarray*}
|a^{(k+1)}_{n+1,n}|&=&|\theta^{(k)}(I_{n-1}-L)(\theta^{(k)})^*+
a^{(k)}_{n+1,n}(\lambda_n-1)\overline{a^{(k)}_{n+1,n+1}}+a^{(k)}_{n+1,n+1}(\lambda_{n+1}-1)\overline{a^{(k)}_{n+1,n}}|\\
&\leq & M_g|a^{(k)}_{n+1,n}a^{(k)}_{n+1,n+1}|.
\end{eqnarray*}

By Lemma \ref{lemcu} and the above two inequalities, we have
\begin{eqnarray}\label{ppi}
|a^{(k+1)}_{n,n+1}a^{(k+1)}_{n+1,n}|\leq
M_g^2|a^{(k)}_{n,n}a^{(k)}_{n+1,n+1}||a^{(k)}_{n,n+1}a^{(k)}_{n+1,n}|\leq
M_g^2(1+|a^{(k)}_{n,n+1}a^{(k)}_{n+1,n}|^{1/2})^2|a^{(k)}_{n,n+1}a^{(k)}_{n+1,n}|.
\end{eqnarray}

Hence
\begin{eqnarray}\label{induction}
|a^{(k+1)}_{n,n+1}a^{(k+1)}_{n+1,n}|^{1/2}\leq
M_g(1+|a^{(k)}_{n,n+1}a^{(k)}_{n+1,n}|^{1/2})|a^{(k)}_{n,n+1}a^{(k)}_{n+1,n}|^{1/2}.
\end{eqnarray}

In what follows, we partition our proof by the following three
claims.

{\bf Claim 1:}  Under the conditions (\ref{rcond}),
$a^{(k)}_{n,n+1}a^{(k)}_{n+1,n}$ converges to zero as $k\to
\infty.$

 Suppose that  $T_1=M_g(1+|a_{n,n+1}a_{n+1,n}|^{1/2})<1$.
Then
\begin{equation}
|a^{(1)}_{n,n+1}a^{(1)}_{n+1,n}|^{1/2}\leq
T_1|a_{n,n+1}a_{n+1,n}|^{1/2}<|a_{n,n+1}a_{n+1,n}|^{1/2}
 \end{equation}
and
\begin{eqnarray*}
|a^{(2)}_{n,n+1}a^{(2)}_{n+1,n}|^{1/2}&\leq&
M_g(1+|a^{(1)}_{n,n+1}a^{(1)}_{n+1,n}|^{1/2})|a^{(1)}_{n,n+1}a^{(1)}_{n+1,n}|^{1/2}\\
&<&T_1|a^{(1)}_{n,n+1}a^{(1)}_{n+1,n}|^{1/2}\leq
T_1^2|a_{n,n+1}a_{n+1,n}|^{1/2}.
 \end{eqnarray*}
We obtain by induction that for $k\geq 1$,
\begin{equation}
|a^{(k+1)}_{n,n+1}a^{(k+1)}_{n+1,n}|^{1/2}<T_1^{k+1}|a_{n,n+1}a_{n+1,n}|^{1/2},
 \end{equation}
which implies
\begin{equation}
\lim_{k\to \infty}|a^{(k)}_{n,n+1}a^{(k)}_{n+1,n}|^{1/2}=0.
 \end{equation}

Suppose that $T_2=M_g(1+|a_{n,n}a_{n+1,n+1}|^{1/2})<1$. Then by
(\ref{ppi}) and  Lemma \ref{lemcu}
\begin{eqnarray*}
|a^{(1)}_{n,n+1}a^{(1)}_{n+1,n}|^{1/2}&\leq &
M_g|a_{n,n}a_{n+1,n+1}|^{1/2}|a_{n,n+1}a_{n+1,n}|^{1/2}\\
&\leq &
|a_{n,n}a_{n+1,n+1}|^{1/2}M_g(1+|a_{n,n}a_{n+1,n+1}|^{1/2})\leq
|a_{n,n}a_{n+1,n+1}|^{1/2}.
 \end{eqnarray*}
Thus $R=M_g(1+|a^{(1)}_{n,n+1}a^{(1)}_{n+1,n}|^{1/2})<T_2<1$. From
(\ref{induction}), we obtain  by induction that
\begin{equation}
|a^{(k+1)}_{n,n+1}a^{(k+1)}_{n+1,n}|^{1/2}<R^{k}|a^{(1)}_{n,n+1}a^{(1)}_{n+1,n}|^{1/2},
 \end{equation}
which also implies
\begin{equation}
\lim_{k\to \infty}|a^{(k)}_{n,n+1}a^{(k)}_{n+1,n}|^{1/2}=0.
 \end{equation}

The proof of claim 1 is complete.

\medskip

{\bf Claim 2:} If  there exists some integer $k$ such that
\begin{equation}\label{dcond}
a^{(k+1)}_{n,n+1}a^{(k+1)}_{n+1,n}=0,
 \end{equation}
 then $\langle h,g\rangle$ is either elementary or not discrete.

 The condition (\ref{dcond}) can be divided into three cases as
 follows:
$$ (i)\ \  a^{(k+1)}_{n,n+1}=0,\ a^{(k+1)}_{n+1,n}\neq 0;\ \  (ii)\ \  a^{(k+1)}_{n,n+1}\neq 0,\ a^{(k+1)}_{n+1,n}=
0;\ \ (iii)\ \   a^{(k+1)}_{n,n+1}=0,\ a^{(k+1)}_{n+1,n}= 0.$$

If the case (i) holds, then  (\ref{bb}), (\ref{gg}) and (\ref{bba})
imply that
$$\beta^{(k+1)}=(0,\cdots,0)^t,\ \  \gamma^{(k+1)}=(0,\cdots,0),\ \
\overline{a^{(k+1)}_{n+1,n+1}}a^{(k+1)}_{n,n}= 1.$$ This implies
that
$$h_{k+1}(q_{0})=h_{k}gh_k^{-1}(q_{0})=q_{0}.$$
Since $h_k$ is loxodromic for $k\geq 1$,  $h_k$ can not swap $q_0$
and $q_{\infty}$. We have $h_{k}(q_{0})=q_{0}$ and by induction,
$$h_1(q_{0})=hgh^{-1}(q_{0})=q_0.$$
If in addition, we have
$$h_1(q_{\infty})=hgh^{-1}(q_{\infty})=q_{\infty}.$$
Then $h$ keeps the set $\{q_0,q_{\infty}\}$ invariant, which implies
that $\langle h,g \rangle$ is elementary. If
$$h_1(q_{\infty})=hgh^{-1}(q_{\infty})\neq q_{\infty},$$ then $g$ and
$h_1$ share  exactly one fixed point so  $\langle h_1,g \rangle$ is
not discrete by \cite[Theorem 3.1]{kam83}. Hence  $\langle h,g
\rangle$ is not discrete.

Similarly, if the case (ii) holds, we can show that $\langle h,g
\rangle$ is either elementary or not discrete.

If the case (iii) holds, as in case (i) and case (ii), we can obtain
$$h_1(q_{0})=hgh^{-1}(q_{0})=q_0,\ h_1(q_{\infty})=hgh^{-1}(q_{\infty})=q_{\infty},$$
which implies  $\langle h,g \rangle$ is  elementary.

The proof of claim 2 is complete.

\medskip

{\bf Claim 3:}  If  for any $k\geq 1$,
\begin{equation}\label{neq0}
a^{(k+1)}_{n,n+1}a^{(k+1)}_{n+1,n}\neq 0,\ \ \mbox{and}\ \
a^{(k+1)}_{n,n+1}a^{(k+1)}_{n+1,n}\to 0, k\to \infty.
 \end{equation}
then  $\langle h,g \rangle$ is  not discrete.

Suppose condition (\ref{neq0}) holds. Then
$|a^{(k+1)}_{n,n+1}a^{(k+1)}_{n+1,n}|$ is a distinct  strictly
decreasing sequence of positive numbers. By Lemma \ref{lemcu}, we
know that $|a^{(k)}_{n,n}a^{(k)}_{n+1,n+1}|\leq
(1+|a^{(k)}_{n,n+1}a^{(k)}_{n+1,n}|^{1/2})^2$ is bounded for all
$k$. It follows from (\ref{bbaa}) that
\begin{eqnarray}|(\beta^{(k)})^*\alpha^{(k)}|^2\leq
4|a^{(k)}_{n,n}a^{(k)}_{n+1,n+1}||a^{(k)}_{n,n+1}a^{(k)}_{n+1,n}|,
\end{eqnarray}
which implies that \begin{equation}(\beta^{(k)})^*\alpha^{(k)}\to 0,
\ \ \mbox{as}\ \ k\to \infty.\end{equation} Observing that
$$\overline{a^{(k)}_{n+1,n+1}}a^{(k)}_{n,n}= 1+(\beta^{(k)})^*\alpha^{(k)}-\overline{a^{(k)}_{n,n+1}}a^{(k)}_{n+1,n},$$
we have
\begin{equation}\overline{a^{(k)}_{n+1,n+1}}a^{(k)}_{n,n}\to 1, \ \ \mbox{as}\ \ k\to
\infty.\end{equation} It follows from Proposition \ref{proquat} that
\begin{eqnarray*}
 |-\gamma^{(k)}L(\theta^{(k)})^*|&\leq& |\gamma^{(k)}||\theta^{(k)}|=
 2\Re(a^{(k)}_{n,n}\overline{a^{(k)}_{n,n+1}})2\Re(a^{(k)}_{n+1,n}\overline{a^{(k)}_{n+1,n+1}})\\
 &\leq &
2(\Re(a^{(k)}_{n,n}\overline{a^{(k)}_{n,n+1}}a^{(k)}_{n+1,n}\overline{a^{(k)}_{n+1,n+1}})
+|a^{(k)}_{n,n}\overline{a^{(k)}_{n,n+1}}a^{(k)}_{n+1,n}\overline{a^{(k)}_{n+1,n+1}}|).
\end{eqnarray*}
 By (\ref{tolam1}), we have
\begin{equation}\label{ann}|a^{(k)}_{n,n}|\to |\lambda_n|, \ \ \mbox{as}\ \ k\to
\infty.\end{equation}

Similarly, by (\ref{tolam2}) we have
\begin{equation}|a^{(k)}_{n+1,n+1}|\to |\lambda_{n+1}|, \ \ \mbox{as}\ \ k\to
\infty.\end{equation}

By (\ref{bbaa}), we have  $(\beta^{(k)})^*\beta^{(k)}\to 0$ or
$(\alpha^{(k)})^*\alpha^{(k)}\to 0$ as $k\to \infty$, which implies
that
\begin{equation}\alpha^{(k)}(\beta^{(k)})^*\to O_{n-1},\ \beta^{(k)}(\alpha^{(k)})^*\to O_{n-1} \ \ \mbox{as}\ \ k\to
\infty,\end{equation} where $O_{n-1}$ is $(n-1)\times (n-1)$ matrix
with all entries equaling $0$.   By (\ref{AA}), we have
\begin{equation}A^{(k)}(A^{(k)})^*\to I_{n-1} \ \ \mbox{as}\ \ k\to
\infty.\end{equation}

Since $|a^{(k+1)}_{n,n+1}|\leq
 M_g|a^{(k)}_{n,n+1}a^{(k)}_{n,n}|$ and $M_g<1$,  by (\ref{ann}) we can choose a positive integer
 $N$ and some $0<s<1$ so that for all $k\geq N$, $|a^{(k+1)}_{n,n+1}|<
 s|\lambda_n||a^{(k)}_{n,n+1}|$, that is
\begin{equation}\frac{|a^{(k+1)}_{n,n+1}|}{|\lambda_n|^{k+1}}<s\frac{|a^{(k)}_{n,n+1}|}{|\lambda_n|^{k}},\ \ \mbox{for all}\ \  k\geq N.\end{equation}
Therefore
\begin{equation}\label{qto1}\lim_{k\to \infty}\frac{|a^{(k)}_{n,n+1}|}{|\lambda_n|^{k}}=0.\end{equation}
Similarly,  we have
\begin{equation}\label{qto2}\lim_{k\to \infty}\frac{|a^{(k)}_{n+1,n}|}{|\lambda_{n+1}|^{k}}=0.\end{equation}

Following J{\o}gensen, we now define the sequence
$f_k=g^{-k}h_{2k}g^{k}$.   As a matrix in $Sp(n,1)$ this is given by
\begin{equation}f_k=\left(
                  \begin{array}{ccc}
                    L^{-k}A^{(2k)}L^k & L^{-k}\alpha^{(2k)}\lambda_n^k & L^{-k}\beta^{(2k)}\lambda_{n+1}^k \\
                    \lambda_n^{-k}\gamma^{(2k)}L^k & \lambda_n^{-k}a^{(2k)}_{n,n}\lambda_n^{k} & \lambda_n^{-k}a^{(2k)}_{n,n+1}\lambda_{n+1}^k \\
                    \lambda_{n+1}^{-k}\theta^{(2k)}L^k& \lambda_{n+1}^{-k}a^{(2k)}_{n+1,n}\lambda_n^{k}
                    & \lambda_{n+1}^{-k}a^{(2k)}_{n+1,n+1}\lambda_{n+1}^k \\
                  \end{array}
                \right).
\end{equation}
Since $|\lambda_{n}\lambda_{n+1}|=1$, by (\ref{qto1}) and
(\ref{qto2}) we have that
\begin{equation}\lim_{k\to \infty}|\lambda_n^{-k}a^{(2k)}_{n,n+1}\lambda_{n+1}^k|=0\ \
\mbox{and}\ \ \lim_{k\to
\infty}|\lambda_{n+1}^{-k}a^{(2k)}_{n+1,n}\lambda_n^{k}|=0\end{equation}
By (\ref{aa}), we have
\begin{equation}|L^{-k}\alpha^{(2k)}\lambda_n^k|^2=\overline{\lambda_n^k}(\alpha^{(2k)})^*(L^{-k})^*L^{-k}\alpha^{(2k)}\lambda_n^k
=2\Re(\frac{\overline{a^{(2k)}_{n+1,n}}}{|\lambda_{n+1}|^{2k}}a^{(2k)}_{n,n}),
\end{equation}
which implies that
\begin{equation}\lim_{k\to \infty}L^{-k}\alpha^{(2k)}\lambda_n^k=(0,\cdots, 0)^t.\end{equation}
Similarly, we have
\begin{equation}\lim_{k\to \infty} L^{-k}\beta^{(2k)}\lambda_{n+1}^k=(0,\cdots, 0)^t,\end{equation}
\begin{equation}\lim_{k\to \infty} \lambda_n^{-k}\gamma^{(2k)}L^k=(0,\cdots, 0),\end{equation}
and
\begin{equation}\lim_{k\to \infty}\lambda_{n+1}^{-k}\theta^{(2k)}L^k=(0,\cdots, 0).\end{equation}
Furthermore, we have  \begin{equation}A^{(2k)}(A^{(2k)})^*\to
I_{n-1}, |\lambda_n^{-k}a^{(2k)}_{n,n}\lambda_n^{k}|\to |\lambda_n|,
\, |\lambda_{n+1}^{-k}a^{(2k)}_{n+1,n+1}\lambda_{n+1}^k|\to
|\lambda_{n+1}|,\ \ \mbox{as}\ \ k\to \infty.\end{equation}

From the above limits, we find that $f_k$ are all distinct and
there exists a subsequence which converges to a element in
$\langle h,g \rangle$.  So the group  $\langle h,g \rangle$ is not
discrete. The proof of claim 3 is complete.

 The proof is complete.

\end{document}